\documentclass[12pt]{amsart}
\usepackage{amssymb}
\usepackage{amsfonts}
\usepackage{amsmath}
\usepackage{graphicx}
\usepackage{xcolor}
\usepackage{mathrsfs}
\usepackage{stmaryrd}
\usepackage{epsfig,color}
\usepackage{blindtext}
\usepackage{enumerate}
\usepackage{hyperref}
\usepackage{url}
\usepackage{bbm}
\usepackage{filecontents}
\usepackage{nicefrac,mathtools}
\usepackage{bm}  
\usepackage[nocompress]{cite}
\DeclareGraphicsExtensions{.pdf,.jpeg,.png}
\usepackage{epstopdf}
\usepackage{cancel} %for editing
\usepackage[normalem]{ulem} %for editing
\usepackage{verbatim} %for comment
\usepackage{enumitem} % for alph* enumerate item
\usepackage{tikz-cd}
\usetikzlibrary{cd}
\usepackage{color}
\usepackage[msc-links, lite]{amsrefs}
\usepackage{geometry}
\usepackage{tikz}
\usetikzlibrary{decorations.markings}
\usetikzlibrary{arrows.meta}

\usepackage{extarrows}
\newtheorem{theorem}{Theorem}[section]
\newtheorem{conjecture}[theorem]{Conjecture}
\newtheorem{proposition}[theorem]{Proposition}
\newtheorem{lemma}[theorem]{Lemma}
\newtheorem{corollary}[theorem]{Corollary}

\theoremstyle{definition}
\newtheorem{definition}[theorem]{Definition}
\newtheorem{example}[theorem]{Example}

\newtheorem{remark}[theorem]{Remark}

\numberwithin{equation}{section}

\numberwithin{equation}{section}

\makeatletter
\let\save@mathaccent\mathaccent
\newcommand*\if@single[3]{%
\setbox0\hbox{${\mathaccent"0362{#1}}^H$}%
\setbox2\hbox{${\mathaccent"0362{\kern0pt#1}}^H$}%
\ifdim\ht0=\ht2 #3\else #2\fi
}

\makeatother

\makeatletter
\newcommand*{\transpose}{%
{\mathpalette\@transpose{}}%
}
\newcommand*{\@transpose}[2]{%
\raisebox{\depth}{$\m@th#1\intercal$}%
}
\makeatother
\usetikzlibrary{hobby}

 \usetikzlibrary{decorations}
\makeatletter
\def\pgfutil@Repeat#1#2{#2\ifnum#1>0
  \expandafter\pgfutil@firstofone\else\expandafter\pgfutil@gobble\fi
  {\expandafter\pgfutil@Repeat\expandafter{\the\numexpr#1-1\relax}{#2}}}
\tikzset{
  dash between/.code args={#1 and #2}{%
    \tikz@addoption{%
      \pgfgetpath\currentpath
      \pgfprocessround{\currentpath}{\currentpath}%
      \pgf@decorate@parsesoftpath{\currentpath}{\currentpath}%
      \pgfmathsetlengthmacro\firstpart{(#1)*\pgf@decorate@totalpathlength}%
      \pgfmathsetlengthmacro\secondpart{(#2-(#1))*\pgf@decorate@totalpathlength}%
      \pgfmathsetlengthmacro\thirdpart{(1-(#2))*\pgf@decorate@totalpathlength}%
      \edef\thirdpart{{\thirdpart}{0pt}}%
      \edef\firstpart{{\firstpart}{0pt}}%
      \pgfmathsetlengthmacro\secondpartlength{\pgfkeysvalueof{/tikz/dash between on}
                                            +(\pgfkeysvalueof{/tikz/dash between off})}%
      \pgfmathtruncatemacro\repetitions{\secondpart/\secondpartlength}%
      \pgfmathsetlengthmacro\secondexpand{\secondpart/\repetitions-\secondpartlength}%
      \edef\secondexpand{\the\dimexpr\pgfkeysvalueof{/tikz/dash between off}+\secondexpand\relax}%
      \edef\secondpart{%
        \pgfutil@Repeat{\the\numexpr\repetitions-1\relax}%
          {{\pgfkeysvalueof{/tikz/dash between on}}{\secondexpand}}%
      }%
      \edef\tikz@temp{\firstpart\secondpart\thirdpart}%
      \expandafter\pgfsetdash\expandafter{\tikz@temp}{+0pt}%
    }
  }
}
\makeatother

\tikzset{
  dash between style/.is choice,
  dash between style/dotted/.style        ={dash between on=\pgflinewidth,dash between off=2pt},
  dash between style/densely dotted/.style={dash between on=\pgflinewidth,dash between off=1pt},
  dash between style/loosely dotted/.style={dash between on=\pgflinewidth,dash between off=4pt},
  dash between style/dashed/.style        ={dash between on=3pt,dash between off=2pt},
  dash between style/loosely dashed/.style={dash between on=3pt,dash between off=6pt},
  dash between style/densely dashed/.style={dash between on=3pt,dash between off=2pt},
  dash between style/no/.style={dash between on=0pt, dash between off=1pt},% dirty
  dash between on/.initial=\pgflinewidth,
  dash between off/.initial=2pt,
  middle dotted line/.style={
    thick,
    dash between=.35 and .65}}

\newcommand\CC{\mathbb{C}}

\DeclareMathOperator{\Pro}{Pro}

\DeclareMathOperator{\homo}{Hom}

\DeclareMathOperator{\et}{\acute{e}t}

\DeclareMathOperator{\Int}{Int}

\DeclareMathOperator{\Gal}{Gal}

\def\holim{\qopname\relax m{holim}}

\title{A homotopical consequence of branched covers}
\author{Runjie Hu}
\newcommand{\Addresses}{{% additional braces for segregating \footnotesize
  \bigskip
  \footnotesize

  Runjie Hu, \textsc{Department of Mathematics, Texas A\&M University, College Sta, TX 77843}\par\nopagebreak
  \textit{E-mail address}, \texttt{ runjie.hu@tamu.edu}
}}
\date{}
\begin{document}

\maketitle

\begin{abstract}
We prove that the profinite completion of a pseudomanifold is Artin-Mazur's \'etale homotopy type construction from its branched covers, which was a statement implicitly conjectured by Sullivan in \cite{Sullivan2009}*{p.~272} around 1970. This is a consequence of the existence of enough $K(\pi,1)$ open dense subspaces in a pseudomanifold.
\end{abstract}

\vspace{0.7cm}

In the study of the underlying homotopy type of a space, one separates rational, prime-local, and finite information by localizations and completions (\cite{Sullivan2009}). The rational homotopy theory is well understood by various algebraic models, such as differential graded algebras, coalgebras, Hopf algebras, Lie algebras, etc. (e.g., \cite{Quillen1969}\cite{sullivan1977}). However, the finite completion of a space which encodes all the finite homotopy information is more complicated and harder to understand (in the literature, this is called the profinite completion).

The (pro-)finite completion of spaces is a homotopical analogue of the corresponding construction for groups. Recall that one way to define the (pro-)finite completion $\widehat{G}$ of a group $G$ is to take the inverse system of finite groups $H$ indexed by all homomorphisms $G\rightarrow H$, with $H$ finite. Taking the inverse limit of this system produces (pro-)finite group $\widehat{G}$. The profinite completion $\widehat{X}$ of a space $X$ is constructed by taking the inverse system of finite spaces $F$ indexed by all maps into such $F$'s (see \cite[Chapter 3]{Artin&Mazur1969}), where a finite space means that each of its homotopy groups is finite and only finitely many are nonzero. There is a way to take the inverse limit of this system of spaces in the homotopy sense to obtain a (pro-)finite space $\widehat{X}$ (see \cite[Definition 3.1]{Sullivan2009}, or \cite{Bousfield&Kan}).

The definition of the profinite completion of spaces is very useful in both algebraic geometry and algebraic topology, but the construction of the inverse system is difficult to control. The following conjecture implicitly brought up by Sullivan in \cite{Sullivan2009}*{p.~272, the last two paragraphs} provides a better understanding of the profinite completion of manifolds.

\begin{conjecture}
There is an analogue of the \'etale
site in algebraic geometry associated to any smooth manifold $M$, defined using branched covers, such that the profinite completion of $M$ is the pro space given by Artin-Mazur's construction applied to this ``\'etale site'' over $M$.
\end{conjecture}

The main motivation of this conjecture is to provide a geometric understanding of the \'etale theory of algebraic varieties. For a compact complex variety $X$, an \'etale morphism over $X$ is a local homeomorphism $f:U\rightarrow X$, and any such map can be completed to an algebraic branched covering. Artin and Mazur formally constructed a (pro-)space $X_{\et}$ from all \'etale morphisms over $X$ and proved that $X_{\et}$ is the profinite completion of $X$ (\cite{Artin&Mazur1969}*{Theorem 12.9}). 

In the geometric setting, the complex variety $X$ is a pseudomanifold and the role of ``\'etale morphisms'' is naturally played by ``geometric'' branched covers (see Section \ref{section branched cover}).

This paper proves a more general result which implies this conjecture. We construct an inverse system by applying Artin-Mazur's construction of the \'etale homotopy type of algebraic varieties (\cite{Artin&Mazur1969}) to branched covers over a pseudomanifold. This is essentially a C\v ech-like construction modified by replacing open subsets with finite covering spaces of open dense sets. The following theorem is our main result (see Theorem \ref{main} and Corollary \ref{main2} for a more rigorous statement).

\begin{theorem}[Theorem \ref{main}]\label{Main1}
The Artin-Mazur's  \'etale homotopy type construction $X_{\et}$ of all branched covers over a pseudomanifold $X$ is the profinite completion of $X$.
\end{theorem}

In \cite{Sullivan2009}*{p.~272}, Sullivan proposed using ``transversality and the representation of Eilenberg-MacLane spaces as large symmetric products
of simple spaces'' to prove a crucial step: ``any finite coefficient cohomology class of a manifold could be killed by passing
to a finite branched cover''. We establish this crucial step by providing a more straightforward geometric construction in Lemma \ref{Key}, where we prove that a pseudomanifold has sufficiently many $K(\pi,1)$ dense open subsets. A statement analogous to Lemma \ref{Key} in the algebraic geometry setting serves as a crucial step in proving Artin-Mazur's result for complex algebraic varieties as well (see \cite{Sullivan1974}*{the paragraph below ``Remarks on the proof''}).

Finally, we use our result to explain a well-known example of Riemann surfaces in Example \ref{example: Riemann surface}.

\subsubsection*{Acknowledgements.} I would like to express my heartfelt thanks to my thesis advisor Dennis Sullivan. I wish to express my gratitude to Irina Bobkova, Mark de Cataldo, Jiahao Hu, John Morgan, Jason Starr, Guozhen Wang, Shmuel Weinberger, Zhouli Xu and Siqing Zhang for their valuable suggestions on this article. I wish to express my gratitude to the Simons Foundation International for its financial support during my research on this topic. I am grateful to the referee who provided many invaluable comments on this article.

\subsection*{Terminology}
We briefly review the following terminology used in this article.
\begin{enumerate}[leftmargin=0.25in]
    \item Let $\mathcal{H}_0$ be the homotopy category of pointed connected CW complexes.
    To every pointed connected site $\mathcal{T}$, Artin-Mazur associate a pro-object $\Pi(\mathcal{T})\in \Pro(\mathcal{H}_0)$ (\cite{Artin&Mazur1969}*{p.~112, ``the Verdier functor''}). 
    \item A CW complex $X$ is called \textbf{$\pi$-finite} if $\pi_0(X)$ is finite, and there exists an $n$ such that, for each $x\in X$, $\pi_i(X,x)=0$ for all $i>n$ and $\pi_i(X,x)$ is finite for all $i$. Let $\mathcal{FH}_0\subset \mathcal{H}_0$ be the full subcategory of $\pi$-finite pointed connected CW complexes. The inclusion functor $\Pro(\mathcal{FH}_0)\rightarrow \Pro(\mathcal{H}_0)$ admits a left adjoint $\widehat{(-)}:\Pro(\mathcal{H}_0)\rightarrow \Pro(\mathcal{FH}_0)$, which is called \textbf{Artin-Mazur's profinite completion} (\cite{Artin&Mazur1969}*{Theorem 3.4}).
    \item For any $\pi$-finite CW complex $F$, the representable functor $[-,F]:\mathcal{H}_0\rightarrow \mathcal{SET}$ takes values in in profinite sets $\mathcal{PRO-FSET}$ (\cite{Sullivan2009}*{Proposition 3.1}). For a pro-object $\{K_i\}\in \Pro(\mathcal{FH}_0)$, Sullivan proves that the functor $\varprojlim_i [-, K_i]:\mathcal{H}_0\rightarrow \mathcal{PRO-FSET}$ exists and satisfies Brown's representability axioms (\cite{Sullivan2009}*{Proposition 3.3}). Hence, it is representable by  a pointed connected CW complex $\holim_i K_i$, which is called \textbf{Sullivan's homotopy inverse limit}.
    \item For a pointed connected CW complex $X$, Sullivan's homotopy inverse limit of Artin-Mazur's profinite completion $\widehat{X}$ for $X$ is called \textbf{Sullivan's profinite completion} of $X$.
\end{enumerate}

\begin{remark}
The following is a modern formulation of the results of this article. 
\begin{itemize}[leftmargin=0.25in]
    \item Let $\mathcal{S}$ be the $\infty$-category of spaces and let $\mathcal{S}^{fc}\subset \mathcal{S}$ be the full subcategory of $\pi$-finite spaces (\cite{Lurie-rational-p-profinite}*{Notation 3.6.1}). Then the restriction functor and the inclusion functor form an adjunction $L_{\mathrm{prof}}:\Pro(\mathcal{S})\rightleftarrows \Pro(\mathcal{S}^{fc}):i$. The functor $L_{\mathrm{prof}}$ is called the \textbf{profinite completion} of pro-spaces.
    \item Let $\mathcal{C}$ be a Grothendieck site, and let $\mathrm{Shv}_{\mathrm{hyp}}(\mathcal{C})$ denote the $\infty$-topos of hypersheaves (\cite{Lurie-HTT}*{paragraph above Lemma 6.5.2.9}). Let $\Pi_{\infty}(\mathrm{Shv}_{\mathrm{hyp}}(\mathcal{C}))\in \Pro(\mathcal{S})$ be the shape of $\mathrm{Shv}_{\mathrm{hyp}}(\mathcal{C})$ (\cite{Lurie-HTT}*{Definition 7.1.6.3}).
    \item Using \cite{Hoyoid-higher-Galois-theory}*{Proposition 5.1} and Lemma \ref{Finiteness}, $\Pi_{\infty}(\mathrm{Shv}_{\mathrm{hyp}}(\mathcal{C}))$ is in the essential image of $\Pro(\mathcal{S}^{fc})$. Our main result (Theorem \ref{Main1}) shows that, for an oriented irreducible pseudomanifold $X$, the shape $\Pi_{\infty}(\mathrm{Shv}_{\mathrm{hyp}}(\mathcal{T}_2(X)))$ of the Grothendieck site of branched covers over $X$ (Definition \ref{Def: topology of branched covers}) is the profinite completion of $X$.
\end{itemize}
\end{remark}

\section{Pseudomanifolds and Branched Covers}\label{section branched cover}

For a PL manifold, one may start from any top cell and bash down walls of top-dimensional faces to reach any other top cell. This geometric property exactly defines a pseudomanifold (\cite{Alexander1920}\cite{Fox1957}\cite{Viro1984}\cite{Berstein&Edmonds1978}\cite{Montesinos1985}\cite{Iori&Piergallini2002}). 

\begin{definition}\label{Def: irreducible pseudomianifold}
A \textbf{pseudomanifold} $X$ of dimension $n$ is a compact polyhedron which admits a finite triangulation satisfying the following conditions:
\begin{enumerate}[leftmargin=0.25in]
    \item  each simplex is a face of some $n$-simplex;
    \item each $(n-1)$-simplex is a face of precisely two $n$-simplices.
\end{enumerate}
A triangulation on $X$ satisfying these properties is called a \textbf{compatible triangulation}. A pseudomanifold $X$ is called \textbf{irreducible} if there is a compatible triangulation on $X$ such that,  for any two $n$-simplices, there exists a path $\gamma$ connecting their barycenters such that $\gamma$ has no intersection with simplices of codimension $2$.
\end{definition}

Geometrically, irreducible pseudomanifolds are special cases of manifolds with singularities, which are also called stratified spaces; see \cite[p.~204]{Weinberger1995StratifiedSpaces}.

\begin{definition}
Let $X$ be an irreducible pseudomanifold, and let $Y$ and $Y'$ be pseudomanifolds.
\begin{enumerate}[leftmargin=0.25in]
    \item An \textbf{orientation} on $Y$ is a choice of orientations on each top-dimensional simplex with respect to a compatible triangulation of $Y$ such that the two orientations on each codimension-one simplex induced from the two adjacent top-dimensional simplices are opposite.
    \item An open subset $U$ of a pseudomanifold $Y$ is \textbf{Zariski open} if $Y-U$ is a subcomplex of codimension at least $2$ with respect to a compatible triangulation on $Y$.
    \item Assume that $X$ and $Y$ are oriented. A \textbf{branched cover} $f:Y\rightarrow X$ is a simplicial map, with respect to some compatible triangulations, whose restriction to each top-dimensional simplex of $Y$ is an orientation-preserving homeomorphism onto a top-dimensional simplex of $X$. There exists a maximal Zariski open subset $U$ of $X$ such that $f\vert_{f^{-1}(U)}$ is a finite covering map. This Zariski open subset $U$ is called the \textbf{unbranched range}, and $X-U$ is called the \textbf{branched locus}. $f^{-1}(U)$ is called the \textbf{unbranched domain}. We say that $f$ is \textbf{connected} if $Y$ is irreducible.
    \item A \textbf{branched morphism} $(f:Y\rightarrow X)\rightarrow (f':Y'\rightarrow X)$ between two branched covers is a branched covering map $g:Y\rightarrow Y'$ over $X$.
\end{enumerate}
\end{definition}

Recall the following lemma of Fox.

\begin{lemma}[Fox's Completing Lemma, \cite{Fox1957}*{p.~245, Existence Theorem}]\label{completinglemma}
Any connected finite covering map $g:D\rightarrow X$ over a Zariski open subset $U$ of an oriented irreducible pseudomanifold $X$ of degree $>1$ can be completed to an irreducible branched cover such that $D$ is the unbranched domain.
\end{lemma}

\begin{proof}
Let $V=X-U$. Triangulate $X$ such that $V$ is a subcomplex. We will inductively construct a branched covering $f:Y\rightarrow X$ extending $g$ such that $D\subset Y$ is the unbranched domain.

Let $V_k$ be the $(n-k)$-skeleton of $V$, where $n$ is the dimension of $X$. Suppose we have constructed $f_{k-1}:Y_{k-1}\rightarrow X-V_{k-1}$ for some $k>1$. 

Let $\{\sigma^{n-k}_i\}$ be the set of all $(n-k)$-simplices of $V$. For each $\sigma^{n-k}_i$, let $L_i$ be its link. Define $Y_k$ to be the union of $Y_{k-1}$ and $\Int(\sigma^{n-k}_{i})\times C(f^{-1})(L_{i})$ for all $i$, where $C(-)$ is the cone construction. Then one can extend $f_{k-1}$ to $f_{k}:Y_{k}\rightarrow X-V_{k}$. 
\end{proof}

The remainder of this paper uses the language of Grothendieck topologies and hypercoverings. We recommend \cite[p.~38]{milneLEC} as an introduction to Grothendieck topologies and \cite{Artin&Mazur1969}*{Chapter 8} for hypercoverings.

\begin{definition}\label{Def: topology of branched covers}
Let $X$ be an irreducible oriented pseudomanifold. Define $\mathcal{T}_2(X)$ to be the following Grothendieck topology. It is the category whose objects are triples $(Y,f,S_f)$ with $f:Y\rightarrow X$ a finite disjoint union of irreducible branched covers and $S_f$ a codimension at least $2$ subpolyderon in $Y$ such that $Y-S_f$ is contained in the unbranched domain of $f$. A morphism $h:(Z,f,S_f)\rightarrow (Y,g,S_g)$ is a branched map $h:Z\rightarrow Y$ over $X$ such that $h(S_f)$ contains $S_g$. A family of morphisms $\{h_i:(Z_i,f_i,S_{f_i})\rightarrow (Y,g,S_{g})\}$ in $\mathcal{T}_2(X)$ is a covering if the family is finite and $\bigcup_ih_i(Z_i-S_{f_i})=Y-S_{g}$.
\end{definition}

To check that $\mathcal{T}_2(X)$ is a Grothendieck topology, one needs to verify that pullbacks exist in $\mathcal{T}_2(X)$. Consider a diagram of branched coverings over $X$ as follows.
\[
\begin{tikzcd}
    & \bigsqcup_{i} (Z'_i,S_{g\circ f'_i}) \arrow[d,"\bigsqcup_i f'_i"] \\
    \bigsqcup_{j} (Z_j,S_{g\circ f_j}) \arrow[r,"\bigsqcup f_j"] & (Y_,S_{g})
\end{tikzcd}
\]
Let $U_{ij}:=(Z'_i-(S_{g\circ f'_i}\cup (f'_i)^{-1}(S_g\cup f_j(S_{g\circ f_j}))))\times_{Y-(S_g\cup f'_i(S_{g\circ f'_i}\cup f_j(S_{g\circ f_j})))} (Z_j-(S_{g\circ f_j}\cup (f_j)^{-1}(S_g\cup f'_i(S_{g\circ f'_i}))))$. Using Lemma \ref{completinglemma}, one completes $U_{ij}$ to a finite disjoint union $W_{ij}$ of irreducible pseudomanifolds together with branched morphisms to $Z'_i$ and $Z_j$ over $X$. Taking the disjoint union $\bigsqcup_{ij}W_{ij}$ gives the pullback of the diagram. 

\section{Key Lemma}

The following lemma is crucial to prove the main result, Theorem \ref{Main1}.

\begin{lemma}[Key Lemma]\label{Key}
Let $x$ be a point of an irreducible oriented pseudomanifold $X$. For any Zariski open neighborhood $V$ of $x$, there exists a smaller Zariski open neighborhood $U$ of $x$ in $V$ such that 

\centering{($\ast$): $U$ is a $K(\pi,1)$ space with $\pi$ a free group.}
\end{lemma}

\begin{proof}
Fix a triangulation of $X$ such that $B=X-V$ is a subcomplex of $X$. By definition, the codimension of $B$ is at least $2$.

It suffices to construct finitely many Zariski open subsets $\{W_i\}$ of $X$ such that each $W_i$ has the property ($\ast$) and $\{W_i\}$ covers $V$. By \cite[Theorem 3.8]{Rourke&Sanderson}, there exists a sequence of regular neighborhoods $N_1\supset N_2\supset ...$ of $B$, whose intersection is $B$, such that the boundary of each $N_i$ is bicollared in $N_{i-1}$ and there exist PL homeomorphisms $N_{i-1}\rightarrow N_i$ which fix $N_{i+1}$. Let $V'$ be the closed complement of a sufficiently small $N_i$. Then it suffices to prove the following claim.

\textbf{Claim}: There exists a finite Zariski open cover $\{W'_j\}$ of $V'$ such that each $W'_j$ satisfies ($\ast$). 

\textbf{Proof of Claim}:
First, fix a compatible triangulation of $X$ such that $V'$ is a subcomplex and $\partial V'$ is bicollared in $X$. 

Choosing a barycenter for each simplex, we obtain a dual decomposition on $V'$. It consists of the cones on the links around the simplices. Note that the complement of the codimensional-two skeleton of the dual decomposition of $V'$ contracts to the $1$-dimensional skeleton. In particular, the complement of the codimensional $2$ skeleton of the dual decomposition of $V'$ is a $K(\pi,1)$-space.

Let $n$ be the dimension of $V'$. We will slightly isotope $V'$ $n$ times such that the union of the images of the complement of the codimensional $2$ skeleton of the dual decomposition is $V'$. That is, the intersection of the images of codimensional $2$ skeletons of the dual decompositions with respect to all compatible triangulations of $X$ is empty.

For each $2$-dimensional simplex $\sigma^2$, one can isotope the interior $n$ times, the intersection of the images of the barycenters is obviously empty.

Inductively assume that we have done these istropies for the $(i-1)$-skeleton of $V'$. Let $\tau^i$ be an arbitrary $i$-dimensional simplex. Let $y_{\tau}$ be its barycenter. Let $A_1,...,A_n$ be the images of the isotopies of the codimensional $2$ skeleton of the dual decomposition in $\partial \tau$. By induction, the intersection of $A_i$'s
is empty. Let $C(A_i)$ be the cone on $A_i$ in $\tau$ with the cone point $y_{\tau}$.

Since $\tau$ is the cone of $\partial \tau$, with the cone point $y_{\tau}$, the $j$-th isotopy of $\partial \tau$ can be extended to an isotopy of $\tau$ such that $C(A_j)$ is the image of $C(A_1)$.

By the general position theorem (see \cite{Rourke&Sanderson}*{p.~61, 5.3 Addendum}), one can inductively slightly isotope $\tau^i$ $n$ times, while fixing a neighborhood of the boundary $\partial \tau^i$, such that the image of each $C(A_j)$ is in general position with the intersection of the images of $C(A_1),...,C(A_{j-1})$. Computing the dimension of the intersection of the images of $A_1,...,A_n$, the intersection must be empty. 
Concatenating the $j$-th isotopy of $\tau$ in this paragraph with the $j$-th isotopy in the previous paragraph produces $n$ small isotopies on the $i$-skeleton of $V'$. This finishes the claim.
\end{proof}

\section{Main Result}

This section proves Theorem \ref{main}, which is Theorem \ref{Main1} in the introduction. Besides the Grothendieck topology $\mathcal{T}_2(X)$, we will consider two other topologies on the irreducible pseudomanifold $X$.
\begin{itemize}
    \item Let $\mathcal{T}_0(X)$ be the usual topology of open subsets.
    \item Let $\mathcal{T}_1(X)$ be the Grothendieck topology of finite covering spaces of open subsets.
\end{itemize}

See \cite{Artin&Mazur1969}*{p.~111} for the definition of a \textit{connected} Grothendieck topology.

\begin{proposition}
The Grothendieck topologies $\mathcal{T}_0(X)$ and $\mathcal{T}_1(X),\mathcal{T}_2(X)$ are all connected.   
\end{proposition}

\begin{proof}
It is immediate that $\mathcal{T}_0(X),\mathcal{T}_1(X)$ are connected. $\mathcal{T}_2(X)$ is connected since we only consider irreducible pseudomanifolds and finite unions of connected branched covers.
\end{proof}

Using the definition of a \textit{homotopy} for two simplicial morphisms $f,g:K_{*}\rightarrow L_*$ between two simplicial objects in a Grothendieck topology given in \cite{Artin&Mazur1969}*{p.~102}, the next result follows directly from \cite{Artin&Mazur1969}*{Corollary 8.13}.

\begin{proposition}
The category $HR(\mathcal{T}_2(X))$ of all hypercoverings of $\mathcal{T}_2(X)$ together with homotopy classes of simplicial morphisms is cofiltering. The same is true for $HR(\mathcal{T}_0(X))$ and $HR(\mathcal{T}_1(X))$.
\end{proposition}

For $i=0,1,2$, let $\pi$ be the connected component functor from $\mathcal{T}_i(X)$ to the category of sets. Then the induced functor on $HR(\mathcal{T}_i(X))$ produces a pro object in the homotopy category of simplicial sets. We denote these pro-spaces by $\{K_0\}$, $\{K_1\}$ and $X_{\et}=\{K_2\}$, respectively.

\begin{definition}\label{Definition etale homotopy}
For an irreducible oriented pseudomanifold $X$, $X_{\et}=\{K_2\}$ is called the \textbf{\'etale homotopy type} of $X$.
\end{definition}

\begin{remark}
This definition mirrors the definition of the \'etale homotopy type of an algebraic variety. Recall that for a proper complex variety $X$, an \'etale morphism is a local homeomorphism $f:U\rightarrow X$. Any such map can be completed to an algebraic branched covering map. The \'etale homotopy type of the variety $X$ is the pro-object in the homotopy category of pointed connected simplicial sets obtained by applying the connected component functor to all hypercoverings of \'etale morphisms over $X$.
\end{remark}

There are continuous maps between these topologies as follows: 
\[
\mathcal{T}_2(X) \xleftarrow{q} \mathcal{T}_1(X) \xrightarrow{q'} \mathcal{T}_0(X)
\]
Here, $q'$ is the inclusion functor from $\mathcal{T}_0(X)$ to $\mathcal{T}_1(X)$, and $q$ is the functor from $\mathcal{T}_2(X)$ to $\mathcal{T}_1(X)$ which maps an object $(Y,f,S_f)$ to $(Y-S_f,f)$. These continuous maps induce natural maps between pro-spaces 
\[
X_{\et}=\{K_2\}\xleftarrow{p} \{K_1\}\xrightarrow{p'} \{K_0\}
\]

\begin{lemma}[\cite{Artin&Mazur1969}*{Theorem 12.1}]\label{Cech nerve homotopy equivalence}
$\{K_0\}$ is homotopy equivalent to $X$.  
\end{lemma}

\begin{lemma}\label{refinement homotopy equivalence}
$p':\{K_1\}\rightarrow \{K_0\}$ is a homotopy equivalence.    
\end{lemma}

\begin{proof}
$\{K_0\}$ can be viewed as a pro subsystem of $\{K_1\}$.
Note that an open covering formed by finite covering spaces of open subsets has some refinement of $\bigsqcup U\rightarrow X$, where each $U\rightarrow X$ is the canonical embedding of an open subset $U$ into $X$. Thus $\{K_0\}$ is cofinal in $\{K_1\}$. This proves that $p'$ is a homotopy equivalence.
\end{proof}

Recall the following definition from \cite{milneLEC}*{p.~48}.

\begin{definition}
A sheaf of sets $\mathcal{F}$ over a Grothendieck topology $G$ is \textbf{locally constant} if there exists a cover $\{U_i\}$ of the terminal object in $G$ such that $\mathcal{F}\vert_{U_i}$ is a constant sheaf for any $i$.
\end{definition}

\begin{lemma}[\cite{Artin&Mazur1969}*{Corollary 10.8}]\label{comparison between singular cohomology and sheaf cohomology}
Each locally constant sheaf of abelian groups $\widetilde{F}$ over $\mathcal{T}_2(X)$ corresponds to a unique local system $\widetilde{F}$ of abelian groups over $\{K_2\}$. Moreover, the usual cohomology $H^*(\{K_2\};\widetilde{F})= \varinjlim H^*(K_2;\widetilde{F})$ is isomorphic to the sheaf cohomology $H^*(\mathcal{T}_2(X);\widetilde{F})$. The same holds for $\mathcal{T}_0(X)$ and $\mathcal{T}_1(X)$.
\end{lemma}

Let $\widetilde{F}$ be a local system of abelian groups over $X$. It induces local systems of abelian groups over $\{K_1\}$ and $\{K_0\}$. It also induces locally constant sheaves over $\mathcal{T}_1(X)$ and $\mathcal{T}_0(X)$. We use the same symbol $\widetilde{F}$ for all these local systems and locally constant sheaves.

\begin{corollary}\label{cohomology and sheaf cohomology}
The sheaf cohomology groups $H^*(\mathcal{T}_1(X);\widetilde{F})$, $H^*(\mathcal{T}_0(X);\widetilde{F})$ and the usual cohomology groups $H^*(\{K_1\};\widetilde{F})= \varinjlim H^*(K_1;\widetilde{F})$, $ H^*(\{K_0\};\widetilde{F})=\varinjlim H^*(K_0;\widetilde{F})$ are all isomorphic.
\end{corollary}

\begin{proof}
This directly follows from Lemma \ref{comparison between singular cohomology and sheaf cohomology}, Lemma \ref{refinement homotopy equivalence} and Lemma \ref{Cech nerve homotopy equivalence}.
\end{proof}

\begin{corollary}\label{sheaf cohomology computes the unbranched domain}
Let $(Y,f,S_f)$ be an object in $\mathcal{T}_2(X)$. Let $q(Y,f,S_f)$ be its image in the category $\mathcal{T}_1(X)$ under $q$. Let $\widetilde{F}$ be a locally constant sheaf over $\mathcal{T}_1(X)$. Then the sheaf cohomology $H^*_{\mathcal{T}_1(X)}(q(Y,f,S_f);\widetilde{F})$ is isomorphic to the usual cohomology $H^*(Y_f-S_f;\widetilde{F})$.
\end{corollary}

\begin{proof}
$q(Y,f,S_f)$ is $f:Y-S_f\rightarrow X$. This corollary follows by applying  Corollary \ref{cohomology and sheaf cohomology} to $Y-S_f$.
\end{proof}

\begin{theorem}[Comparison of Cohomologies]\label{cohomologycomparison}
For any local system of finite groups $\widetilde{F}$ over an irreducible oriented pseudomanifold $X$, there is a canonical isomorphism $H^i(X;\widetilde{F})\cong \varinjlim H^i(K_2;p_* \widetilde{F})$ induced by the map $X\simeq \{K_1\}\xrightarrow{p} \{K_2\}$, where $i=1$ if $F$ is a nonabelian group and $i\geq 1$ if $F$ is an abelian group.
\end{theorem}

\begin{lemma}\label{principalbundlelemma}
For any finite group $G$, the canonical map $X\simeq\{K_1\}\xrightarrow{p} \{K_2\}$ induces a bijection $p^*: \homo(\pi_1(\{K_2\}),G)\rightarrow \homo(\pi_1(X),G)$.
\end{lemma}

\begin{proof}
By \cite{Artin&Mazur1969}*{Corollary 10.7}, $\homo(\pi_1(\{K_2\}),G)$ is bijective to locally constant sheaves of $G$-sets whose stalks are equivariantly bijective to $G$ over $\mathcal{T}_2(X)$ and the same is true for $\{K_1\}$ and $\mathcal{T}_1(X)$. Since $\{K_1\}\simeq \{K_0\}\simeq X$, it suffices to show that locally constant sheaves of finite sets over $\mathcal{T}_2(X)$ are bijective to those over $\mathcal{T}_1(X)$.

Let $\mathcal{F}$ be a locally constant sheaf over $\mathcal{T}_2(X)$. By definition, there exists a cover $(Y_i,f_i,S_{f_i})$ over $(X,Id_X,\emptyset)$ in $\mathcal{T}_2(X)$ such that $\mathcal{F}\vert_{(Y_i,f_i,S_{f_i})}$ is a constant sheaf for any $i$. This means that $q^*\mathcal{F}\vert_{Y_i-S_{f_i}}$ is a constant sheaf on the overcategory $\mathcal{T}_1(X)_{(Y_i-S_{f_i},f_i)}$. Since the union of $f_i(Y_i-S_{f_i})$ is $X$, this shows that $q^*\mathcal{F}$ is a locally constant sheaf over $\mathcal{T}_1(X)$.

Conversely, a locally constant sheaf $\mathcal{G}$ of finite sets over $\mathcal{T}_1(X)$ is the same as a locally constant sheaf $\mathcal{G}'$ of finite sets over $X$ with the usual topology. This means that there exists a finite covering space $f:\widetilde{X}\rightarrow X$ such that $f^*\mathcal{G}'$ is a constant sheaf over $\widetilde{X}$. As $(\widetilde{X},f,\emptyset)$ is an object in $\mathcal{T}_2(X)$, $q_*\mathcal{G}\vert_{(\widetilde{X},f,\emptyset)}$ is a constant sheaf. This completes the proof.
\end{proof}

\begin{proof}[Proof of Theorem \ref{cohomologycomparison}]
Lemma \ref{principalbundlelemma} proves the case of $H^1$ for both nonabelian and abelian sheaves.

For higher degrees, since $X\simeq \{K_0\}\simeq \{K_1\}$, $H^*(X;\widetilde{F})\cong\varinjlim H^*(K_1;\widetilde{F})$. By Corollary \ref{cohomology and sheaf cohomology}, $\varinjlim H^*(K_1;\widetilde{F})$ is isomorphic to the sheaf cohomology $H^*(\mathcal{T}_1(X);\widetilde{F})$.

Consider the Leray spectral sequence for $q: \mathcal{T}_1(X)\rightarrow \mathcal{T}_2(X)$ and the sheaf $\widetilde{F}$: 
\[
H^r(\mathcal{T}_2(X);R^s q_* \widetilde{F})\Rightarrow H^{r+s}(\mathcal{T}_1(X);\widetilde{F})
\] 

It suffices to show that $R^s q_* \widetilde{F}=0$ for all $s>0$. By Corollary \ref{sheaf cohomology computes the unbranched domain} and \cite{milneLEC}*{p.~132, the paragraph above ``Loosely speaking, ...''} this is equivalent to proving the following claim.

\textbf{Claim}: For any object $(Y,f,S_f)$ in $\mathcal{T}_2(X)$, any cohomology class $t$ in the usual cohomology $H^s(Y-S_f;\widetilde{F})$ and any $y$ in $Y-S_f$, there exists a morphism $h:(Z,g,S_g)\rightarrow (Y,f,S_f)$ in $\mathcal{T}_2(X)$ such that $y\in h(Z-S_g)$ and the restriction $h^*t$ of $t$ along $h$ in the usual cohomology $H^s(Z-S_g;\widetilde{F})$ vanishes.

\textbf{Proof of Claim}:
By Lemma \ref{Key}, there exists a $K(\pi,1)$ Zariski open neighborhood $V$ of $y$ in $Y_f-S_f$ with $\pi$ a free group. After passing to some finite covering space $W$ of $V$, we may assume that $\widetilde{F}$ is a constant sheaf over $W$ with the usual topology. Note that $W$ is also a $K(\pi,1)$ space. 

For $s=1$, the restriction of $t$ to $W$ corresponds to a principal bundle with group $F$ over $W$. After replacing $W$ by some finite covering space, the principal bundle is trivial, and hence the restriction of $t$ vanishes.

For $s>1$, since $W$ is homotopy equivalent to a wedge of circles, its higher cohomology groups vanish. Thus, the restriction of $t$ to $W$ vanishes.

Use Lemma \ref{completinglemma} to complete $W\rightarrow Y\rightarrow X$ to a finite union of connected branched covers $g:Z\rightarrow X$ and let $S_g=Z-W$. This finishes the proof of the claim and the theorem.
\end{proof} 

\begin{remark}
The claim in this proof is the sentence ``any finite coefficient cohomology class of a manifold could be killed by passing
to a finite branched cover'' in \cite{Sullivan2009}*{p.~272}.
\end{remark}

Recall from \cite{Artin&Mazur1969}*{Example 9.9} that, for a group $G$, a $G$-set is \textbf{connected} if the $G$-action is transitive. This is equivalent to saying that the $G$-set is isomorphic to $G/H$ for some subgroup $H$ of $G$. The category of $G$-sets, equipped with the relevant coverings, is a connected Grothendieck topology. 
 
\begin{lemma}\label{Finiteness}
Each $\pi_i(\{K_2\})$ is profinite.
\end{lemma}

\begin{proof}
It suffices to prove that each $\pi_i(K_2)$ is finite for any $K_2$. Let $L_*$ be the hypercovering in $\mathcal{T}_2(X)$ which corresponds to the simplicial set $K_2$. For each fixed $i$, we may replace $L_*$ by its $(i+1)$-skeleton without changing $\pi_i(K_2)$. So  we may assume that $L_*$ is isomorphic to its $(i+1)$-skeleton. Then only finitely many connected branched covers $(Y,f,S_f)$ over $X$ appear in $L_*$. 

Let $x_0$ be a basepoint in the intersection $W=\bigcap f(Y-S_f)\subset X$. Then each set $f^{-1}(x_0)$ is finite and has a natural $\pi_1(W,x_0)$-action. The sets $f^{-1}(x_0)$ for $(Y,f,S_f)$ appearing in $L_*$ assemble into a simplicial $\pi_1(W,x_0)$-set $L'_*$. Since only finitely many connected $\pi_1(W,x_0)$-sets appear in $L'_*$, the action of $\pi_1(W,x_0)$ factors through some finite quotient group $G$. Since $L_*$ is a hypercovering in $\mathcal{T}_2(X)$, $L'_*$ is a hypercovering of $G$-sets. 

Let $K'_*$ be the simplicial set obtained by applying the connected component functor levelwise to $L'_*$. Since $K'_*=K_2$ as simplicial sets, the lemma follows from \cite{Artin&Mazur1969}*{Proposition 11.3}.
\end{proof}

\begin{theorem}\label{main}
Let $X$ be an irreducible oriented pseudomanifold. The \'etale homotopy type $X_{\et}=\{K_2\}$ (Definition \ref{Definition etale homotopy}) of the Grothendieck topology $\mathcal{T}_2(X)$ of branched covers over $X$ (Definition \ref{Def: topology of branched covers}) is Artin-Mazur's profinite completion of $X$.
\end{theorem}

\begin{proof}
By Lemma \ref{Finiteness} and \cite{Artin&Mazur1969}*{Theorem 4.3}, it suffices to show that the map $X\simeq \{K_0\}\simeq \{K_1\}\xrightarrow{p}\{K_2\}$ induces isomorphisms $\pi_1(X)^{\wedge}\rightarrow \pi_1(\{K_2\})$ and $H^*(\{K_2\};\widetilde{F})=\varinjlim H^*(K_2;\widetilde{F})\rightarrow H^*(X;\widetilde{F})$ for any local system of finite abelian groups $\widetilde{F}$ over $X$. These two isomorphisms are provided by Theorem \ref{cohomologycomparison} and Lemma \ref{principalbundlelemma}.
\end{proof}

\begin{corollary}\label{main2}
Sullivan's homotopy inverse limit of $\{K_2\}$ is homotopy equivalent to Sullivan's profinite completion (\cite{Sullivan2009}*{Chapter 3}) of $X$.
\end{corollary}

\begin{proof}
Let $X'$ be Sullivan's homotopy inverse limit of $\{K_2\}$. Let $\widehat{X}$ be the profinite completion of $X$, which is a pro space constructed in \cite{Artin&Mazur1969}*{Chapter 3}, and let $\widehat{X}_S$ be the homotopy inverse limit of $\widehat{X}$. Then Sullivan's profinite completion of $X$ is $\widehat{X}_S$. The homotopy equivalence $\widehat{X}\simeq \{K_2\}$ in Theorem \ref{main} induces isomorphisms of pro-groups $\pi_i(\widehat{X})\cong \pi_i(\{K_2\})$ and a map $X'\rightarrow \widehat{X}_S$. By Lemma \ref{Finiteness}, $\pi_i(X')$ is the inverse limit of $\pi_i(\{K_2\})$. Since $\pi_i(\widehat{X}_S)$ is the inverse limit of $\pi_i(\widehat{X})$, this corollary follows.
\end{proof}

We apply this result to a well-known example involving Riemann surfaces.

\begin{example}\label{example: Riemann surface}
Let $X$ be a closed Riemann surface. Then $X$ is biholomorphic to an algebraic curve defined by some complex homogenous polynomials $f_1,...,f_r$ in the complex projective space. Let $\sigma$ be a field automorphism of $\CC$. Let $f^{\sigma}_1,...,f^{\sigma}_r$ be the polynomials obtained by applying $\sigma$ to the coefficients of $f_1,...,f_r$. Let $X^{\sigma}$ be the Riemann surface defined by $f^{\sigma}_1,...,f^{\sigma}_r$.

Any branched cover $\widetilde{X}$ over $X$ is also an algebraic curve. Moreover, the covering map $\phi:\widetilde{X}\rightarrow X$ is algebraic. Let $\mathcal{T}_2(X)$ be Grothendieck topology of branched covers over $X$ as in Definition \ref{Def: topology of branched covers}. Then $\sigma$ induces a homeomorphism of Grothendieck topology $\mathcal{T}_2(X^{\sigma})\rightarrow \mathcal{T}_2(X)$ by mapping $\phi:\widetilde{X}\rightarrow X$ to $\phi^{\sigma}:\widetilde{X}^{\sigma}\rightarrow X^{\sigma}$. \qed
\end{example}

As a corollary of this example, one obtains the following known result.

\begin{corollary}
Every Galois automorphism of $\CC$ preserves the profinite homotopy type of a closed Riemann surface. In particular, if $X$ is a Riemann surface defined over a subfield $K\subset \CC$, then there is a canonical homomorphism from $\Gal(\CC/K)$ to the group of homotopy classes of self homotopy equivalences of the profinite completion of $X$.
\end{corollary}

\bibliographystyle{amsalpha}
\bibliography{ref}

\Addresses

\end{document}